\def\spacingset#1{\renewcommand{\baselinestretch}{#1}\small\normalsize}

\documentstyle[12pt]{article}
\input mssymb.tex

\begin{document}
% Define all the custom commands here.
% === LaTeX commands
\newcommand{\nc}{\newcommand}
\nc{\rnc}{\renewcommand}
\nc{\nt}{\newtheorem}
\nc{\be}{\begin}
\nc{\erf}[1]{$\ (\ref{#1}) $}
\nc{\rf}[1]{$\ \ref{#1} $}
\nc{\lb}[1]{\mbox {$\label{#1}$} }
\nc{\hr}{\hrulefill}
\nc{\noi}{\noindent}

\nc{\eq}{\begin{equation}}
\nc{\en}{\end{equation}}
\nc{\eqa}{\begin{eqnarray}}
\nc{\ena}{\end{eqnarray}}

\nc{\ra}{\rightarrow}
\nc{\la}{\leftarrow}
\nc{\da}{\downarrow}
\nc{\ua}{\uparrow}
\nc{\Ra}{\Rightarrow}
\nc{\La}{\Leftarrow}
\nc{\Da}{\Downarrow}
\nc{\Ua}{\Uparrow}

\nc{\uda}{\updownarrow}
\nc{\Uda}{\Updownarrow}
\nc{\lra}{\longrightarrow}
\nc{\lla}{\longleftarrow}
\nc{\llra}{\longleftrightarrow}
\nc{\Lra}{\Longrightarrow}
\nc{\Lla}{\Longleftarrow}
\nc{\Llra}{\Longleftrightarrow}

\nc{\mt}{\mapsto}
\nc{\lmt}{\longmapsto}
\nc{\lt}{\leadsto}
\nc{\hla}{\hookleftarrow}
\nc{\hra}{\hookrightarrow}
\nc{\lgl}{\langle}
\nc{\rgl}{\rangle}

% === "Graphic" commands
\nc{\stla}{\stackrel{d}{\la}}
\nc{\pard}{\partial \da}
\nc{\gdot}{\circle*{0.5}}

% === page layout commands

\rnc{\baselinestretch}{1.2}      % "double" space
\nc{\bl}{\vspace{1ex}}           % leave a blank line
% === counter commands
%\rnc{\theequation}{\arabic{chapter}.\arabic{equation}}  % change eqn #
\rnc{\theequation}{\arabic{section}.\arabic{equation}}  % change eqn #
\newcounter{xs}
\newcounter{ys}
\newcounter{os}
\nt{thm}{Theorem}[section]
\nt{dfn}[thm]{Definition}
\nt{pro}[thm]{Proposition}
\nt{cor}[thm]{Corollary}
\nt{con}[thm]{Conjecture}
\nt{lem}[thm]{Lemma}
\nt{rem}[thm]{Remark}
% === text abbreviations
\nc{\Poincare}{\mbox {Poincar$\acute{\rm e}$} }

% === mathematical abbreviations
\nc{\bA}{\mbox{${\bf A}$\ }}
\nc{\bB}{\mbox{${\bf B}$\ }}
\nc{\bC}{\mbox{${\bf C}$\ }}
\nc{\bD}{\mbox{${\bf D}$\ }}
\nc{\bE}{\mbox{${\bf E}$\ }}
\nc{\bF}{\mbox{${\bf F}$\ }}
\nc{\bG}{\mbox{${\bf G}$\ }}
\nc{\bH}{\mbox{${\bf H}$\ }}
\nc{\bI}{\mbox{${\bf I}$\ }}
\nc{\bJ}{\mbox{${\bf J}$\ }}
\nc{\bK}{\mbox{${\bf K}$\ }}
\nc{\bL}{\mbox{${\bf L}$\ }}
\nc{\bM}{\mbox{${\bf M}$\ }}
\nc{\bN}{\mbox{${\bf N}$\ }}
\nc{\bO}{\mbox{${\bf O}$\ }}
\nc{\bP}{\mbox{${\bf P}$\ }}
\nc{\bQ}{\mbox{${\bf Q}$\ }}
\nc{\bR}{\mbox{${\bf R}$\ }}
\nc{\bS}{\mbox{${\bf S}$\ }}
\nc{\bT}{\mbox{${\bf T}$\ }}
\nc{\bU}{\mbox{${\bf U}$\ }}
\nc{\bV}{\mbox{${\bf V}$\ }}
\nc{\bW}{\mbox{${\bf W}$\ }}
\nc{\bX}{\mbox{${\bf X}$\ }}
\nc{\bY}{\mbox{${\bf Y}$\ }}
\nc{\bZ}{\mbox{${\bf Z}$\ }}
\nc{\cA}{\mbox{${\cal A}$\ }}
\nc{\cB}{\mbox{${\cal B}$\ }}
\nc{\cC}{\mbox{${\cal C}$\ }}
\nc{\cD}{\mbox{${\cal D}$\ }}
\nc{\cE}{\mbox{${\cal E}$\ }}
\nc{\cF}{\mbox{${\cal F}$\ }}
\nc{\cG}{\mbox{${\cal G}$\ }}
\nc{\cH}{\mbox{${\cal H}$\ }}
\nc{\cI}{\mbox{${\cal I}$\ }}
\nc{\cJ}{\mbox{${\cal J}$\ }}
\nc{\cK}{\mbox{${\cal K}$\ }}
\nc{\cL}{\mbox{${\cal L}$\ }}
\nc{\cM}{\mbox{${\cal M}$\ }}
\nc{\cN}{\mbox{${\cal N}$\ }}
\nc{\cO}{\mbox{${\cal O}$\ }}
\nc{\cP}{\mbox{${\cal P}$\ }}
\nc{\cQ}{\mbox{${\cal Q}$\ }}
\nc{\cR}{\mbox{${\cal R}$\ }}
\nc{\cS}{\mbox{${\cal S}$\ }}
\nc{\cT}{\mbox{${\cal T}$\ }}
\nc{\cU}{\mbox{${\cal U}$\ }}
\nc{\cV}{\mbox{${\cal V}$\ }}
\nc{\cW}{\mbox{${\cal W}$\ }}
\nc{\cX}{\mbox{${\cal X}$\ }}
\nc{\cY}{\mbox{${\cal Y}$\ }}
\nc{\cZ}{\mbox{${\cal Z}$\ }}

\nc{\rightcross}{\searrow \hspace{-1 em} \nearrow}
\nc{\leftcross}{\swarrow \hspace{-1 em} \nwarrow}
\nc{\upcross}{\nearrow \hspace{-1 em} \nwarrow}
\nc{\downcross}{\searrow \hspace{-1 em} \swarrow}
\nc{\prop}{| \hspace{-.5 em} \times}
\nc{\wh}{\widehat}
\nc{\wt}{\widetilde}
\nc{\nonum}{\nonumber}
\nc{\nnb}{\nonumber}
 \nc{\half}{\mbox{$\frac{1}{2}$}}
\nc{\Cast}{\mbox{$C_{\frac{\infty}{2}+\ast}$}}
\nc{\Casth}{\mbox{$C_{\frac{\infty}{2}+\ast+\frac{1}{2}}$}}
\nc{\Casm}{\mbox{$C_{\frac{\infty}{2}-\ast}$}}
\nc{\Cr}{\mbox{$C_{\frac{\infty}{2}+r}$}}
\nc{\CN}{\mbox{$C_{\frac{\infty}{2}+N}$}}
\nc{\Cn}{\mbox{$C_{\frac{\infty}{2}+n}$}}
\nc{\Cmn}{\mbox{$C_{\frac{\infty}{2}-n}$}}
\nc{\Ci}{\mbox{$C_{\frac{\infty}{2}}$}}
\nc{\Hast}{\mbox{$H_{\frac{\infty}{2}+\ast}$}}
\nc{\Hasth}{\mbox{$H_{\frac{\infty}{2}+\ast+\frac{1}{2}}$}}
\nc{\Hasm}{\mbox{$H_{\frac{\infty}{2}-\ast}$}}
\nc{\Hr}{\mbox{$H_{\frac{\infty}{2}+r}$}}
\nc{\Hn}{\mbox{$H_{\frac{\infty}{2}+n}$}}
\nc{\Hmn}{\mbox{$H_{\frac{\infty}{2}-n}$}}
\nc{\HN}{\mbox{$H_{\frac{\infty}{2}+N}$}}
\nc{\HmN}{\mbox{$H_{\frac{\infty}{2}-N}$}}
\nc{\Hi}{\mbox{$H_{\frac{\infty}{2}}$}}
\nc{\Ogast}{\mbox{$\Omega_{\frac{\infty}{2}+\ast}$}}
\nc{\Ogi}{\mbox{$\Omega_{\frac{\infty}{2}}$}}
\nc{\Wedast}{\bigwedge_{\frac{\infty}{2}+\ast}}

\nc{\Fen}{\mbox{$F_{\xi,\eta}$}}
\nc{\Femn}{\mbox{$F_{\xi,-\eta}$}}
\nc{\Fp}{\mbox{$F_{0,p}$}}
\nc{\Fenp}{\mbox{$F_{\xi',\eta'}$}}
\nc{\Fuv}{\mbox{$F_{\mu,\nu}$}}
\nc{\Fuvp}{\mbox{$F_{\mu',\nu'}$}}
\nc{\cpq}{\mbox{$c_{p,q}$}}
\nc{\Drs}{\mbox{$\Delta_{r,s}$}}
\nc{\spq}{\mbox{$\sqrt{2pq}$}}
\nc{\Mcd}{\mbox{$M(c,\Delta)$}}
\nc{\Lcd}{\mbox{$L(c,\Delta)$}}
\nc{\Wxv}{\mbox{$W_{\chi,\nu}$}}
\nc{\vxv}{\mbox{$v_{\chi,\nu}$}}
\nc{\dd}{\mbox{$\widetilde{D}$}}

\nc{\diff}{\mbox{$\frac{d}{dz}$}}
\nc{\Lder}{\mbox{$L_{-1}$}}
\nc{\bone}{\mbox{${\bf 1}$}}
\nc{\px}{\mbox{${\partial_x}$}}
\nc{\py}{\mbox{${\partial_y}$}}
 % === the main sections
\setcounter{equation}{0}

\vspace{1in}
\begin{center}
{{\LARGE\bf
Commutative Quantum Operator Algebras
 }}
\end{center}
\addtocounter{footnote}{0}
\footnotetext{1991 Mathematics Subject Classification. Primary 81T70, 17B68.}
% 81T70 is Quantization in field theory; cohomological methods.
% 17B68 is Virasoro and related algebras.
\vspace{1ex}
\begin{center}
Bong H. Lian and Gregg J. Zuckerman
\end{center}
\addtocounter{footnote}{0}
\footnotetext{B.H.L. is supported by grant DE-FG02-88-ER-25065.
G.J.Z. is supported by
NSF Grant DMS-9307086 and DOE Grant DE-FG02-92-ER-25121.}

% === ABSTRACT
\vspace{1ex}
\begin{quote}
{\footnotesize
ABSTRACT. A key notion bridging the gap between {\it quantum operator algebras}
\cite{LZ10} and {\it vertex operator algebras} \cite{Bor}\cite{FLM} is the
definition of the commutativity of a pair of quantum operators
(see section 2 below). This is not commutativity in any ordinary sense, but it
is clearly the correct generalization to the quantum context. The main purpose
of the current paper is to begin laying the foundations for a complete
mathematical theory of {\it commutative quantum operator algebras.}  We give
proofs of most of the relevant results announced in \cite{LZ10}, and we carry
out some calculations with sufficient detail to enable the interested reader to
become proficient with the algebra of commuting quantum operators.
}
\end{quote}
\begin{quote}
{\it We dedicate this paper to the memory of Feza G\"ursey.}
\end{quote}
\addtocounter{footnote}{0}

\section{ Introduction}

\ \ \ \ \ Formal infinite sums of linear operators in a vector space have been
explicitly appearing in quantum field theory since the late 1920's.  A
particular class of such formal sums is studied by the authors \cite{LZ10} in
the context of conformal field theory and string theory.  Let $V$ be a $\bZ$
 doubly graded vector space $V=\oplus V^n[m]$.  Let $z$ be a formal variable
(later to be thought of as a point in the punctured complex plane).  A
homogeneous quantum operator is a formal power series $u(z)=\sum u(n)z^{-n-1}$
of degrees, say $d_1,d_2$, where the coefficients $u(n)$ are linear maps in $V$
with degrees $d_1,d_2-n-1$.
A quantum operator is a finite sum of homogeneous quantum operators, and we
denote the space of all quantum operators in $V$ by $QO(V)$.

The famous vertex operator construction of the old dual resonance model
provides the first concrete and nontrivial examples of quantum operators to
appear in the physics literature, and later in the mathematics literature. Here
the space $V$ is the Fock space, which is only singly graded by weight. The
famous construction of Virasoro yields a quantum operator $L(z)$ acting in $V$
and having the property that the coefficients $L(n)$ of $L(z)$ span, together
with the identity operator, a central extension of the Lie algebra of Fourier
polynomial vector fields on the circle.

The BRST quantization of classical strings provides examples of quantum
operators with nonzero ghost number.  Now the space $V$ is a tensor product of
a so-called ghost Fock space with the state space of a conformal field theory.
A special quantum operator known as the BRST current $J(z)$ plays a key role:
the coefficient $J(0)$ becomes, under suitable conditions, a cohomology
operator in the space $V$.  The much studied BRST cohomology groups
\cite{KO}\cite{FMS}\cite{Fe}\cite{FGZ}
\cite{LZ3}\cite{LZ4}\cite{LZ5} are then defined to be the cohomology groups
associated to the differential, $J(0)$. Next to $J(z)$ the most important
quantum operator in the theory is the ghost field $b(z)$, whose coefficient
$b(1)$ induces in the cohomology the structure of a Batalin-Vilkovisky algebra.
Much of the recent interest in BRST string theory has centered on the myriad
algebraic structures possessed both by BRST cohomology as well as the complex
$V$ itself.

Both the dual resonance theory and the BRST-quantized string theory deal with
the above-mentioned quantum operators as distinguished members of certain
infinite dimensional linear systems of quantum operators. The theory of vertex
operator algebras (VOAs) \cite{Bor}\cite{FLM} establishes a remarkable but
rather complex foundation for the study of such systems. In \cite{LZ9}, the
authors successfully adapted the theory of VOAs to the theory of BRST
cohomology.  In an attempt to better understand this application of VOA theory,
the authors recently introduced a more general class of systems called quantum
operator algebras (QOAs) \cite{LZ10}.
  The latter can also be viewed as elementary generalizations of the operator
algebras well-known in linear algebra and mathematical physics. We want to
emphasize that the operators, rather than the states, are fundamental in our
new point of view.

A key notion bridging the gap between QOAs and VOAs is the definition of
commutativity of a pair of quantum operators \cite{LZ10} (see also
\cite{FLM}\cite{FHL}\cite{dl}\cite{li}). This is not commutativity in any
ordinary sense, but it is clearly the correct generalization to the quantum
context. The main purpose of the current paper is to begin laying the
foundations for a complete mathematical theory of commutative quantum operator
algebras.  We give proofs of most of the relevant results announced in
\cite{LZ10}, and we carry out some calculations with sufficient detail to
enable the interested reader to become proficient with the algebra of commuting
quantum operators.

The main features of this operator calculus were discovered by the quantum
physicists.  These features are certainly captured by the theory of VOAs, and
we have benefitted enormously from the insights in the VOA literature,
especially with regard to the notion of commutativity. However, because VOA
theory emphasizes states in $V$ over quantum operators in $QO(V)$, we feel that
much of the beauty and simplicity of the physicists' operator calculus gets
lost in the translation.  Moreover, we have found that the construction of the
main nontrivial examples of VOAs becomes much simpler and much closer to the
original physical inspiration when carried out in the language of commutative
quantum operator algebras.  Thus, a subsidiary purpose of this paper is to
present in fairly great detail the construction of the simplest CQOAs that
enter into the BRST construction.

Here is a brief summary of the contents of this paper:

In section 2, we state the definitions of our main concepts:  quantum
operators, matrix elements, Wick products, iterated Wick products, the
infinitely many ``circle'' products, the operator product expansion, the
notions of locality and commutativity, and finally the corresponding notions of
local and commutative quantum operator algebras. We discuss some elementary
known facts about commutativity.  We then define the notion of a semi-infinite
commutative
algebra, abstracted from our theory of CQOAs.

In section 3, we introduce what we call the Wick calculus, which deals with
operator products of the form $:t(z)u(z):v(w)$ as well as $t(z):u(w)v(w):$
under the assumption that the quantum operators $t(z)$, $u(z)$, and $v(z)$ are
{\it pairwise commutative}. Here, $:t(z)u(z):$ denotes the Wick or normal
ordered product of $t(z)$ with $u(z)$.  The Wick calculus is essential for both
computations as well as for theoretical issues, such as the explicit
construction of CQOAs.  Section 3 continues with the construction of the CQOA
$O(b, c)$, which acts in the ghost Fock space of the BRST construction.  This
section concludes with a construction of the CQOA $O_\kappa(L)$, which arose
originally in the seminal work of BPZ \cite{BPZ}, and which acts in the state
space of any conformal field theory having central charge $\kappa$. Both
examples are special, in that the two algebras are spanned by iterated Wick
products of derivatives of finitely many generating quantum operators. In fact,
we exhibit explicit bases consisting of such iterated products.

In section 4, we discuss the BRST construction in the language of what we call
conformal QOAs.  Given a conformal QOA $O$ with central charge $\kappa$,
we form the tensor product $C^*(O) = O(b, c)\otimes O$.
We then construct the special quantum operator $J(z)$, referred to above as the
BRST current. We also give a simple characterization of $J(z)$.  We then recall
the famous result that the coefficient $J(0)$ is square-zero if and only if
$\kappa = 26$.
After that we specialize to this central charge.  As in our past paper
\cite{LZ10}, there is a change in point of view:  rather than discussing $J(0)$
as a differential in the BRST state space we instead consider $[J(0),-]$ as a
square-zero derivation of the BRST quantum operator algebra, $C^*(O)$.  Thus,
we view BRST quantization in the operator picture, which is closer to the
original physical context.

A particular feature of section 4 is our calculation for arbitrary central
charge $\kappa$ of the graded commutator $[J(0), J(0)]= 2J(0)^2$, which can be
obtained from the calculation of the operator product, $J(z)J(w)$. Appealing to
the abstract development in section 3, we carry out the evaluation of this
special operator product via the Wick calculus.  This computation provides a
model for the kinds of calculations made routinely by conformal field theorists
since the pioneering work of BPZ \cite{BPZ}.

The main result of section 4 is Theorem \ref{5.3}, which states that the Wick
product induces a graded commutative associative product on the cohomology of
$C^*(O)$ with respect to the derivation, $[J(0),-]$.  This theorem first
appeared in work of E. Witten \cite{Wi3}, who called the ghost number zero
subalgebra the ``ground ring of a string background''.  An approach to this
theorem via VOA theory appears in \cite{LZ9}.  The approach in the current
paper is via CQOA theory, and hopefully appears as a significant simplication
of our earlier work, as well as a clear exposition of Witten's work.

In section 5, we develop the theory of the ghost field, $b(z)$, and its
coefficient, $b(1)$.  As a preparation, we remind the reader of the definition
of a Batalin-Vilkovisky (BV) operator and BV algebra. The main result of
section 5 is Theorem \ref{5.4}, which states that the operator $b(1)$ induces a
BV operator acting in the BRST cohomology algebra. Thus the cohomology becomes
a BV algebra.  This theorem was inspired by work of Witten and Zwiebach
\cite{WZ}, and first appeared in \cite{LZ9}, where it was derived via
identities from VOA theory.  Again, we hope that our new approach via CQOA
theory sheds light on both our original discovery as well as the physical
inspiration.

 Nontrivial examples of BV structure in BRST cohomology appear in section 3 of
\cite{LZ9} as well as section 6 of \cite{LZ10}. We plan to return to these
explicit examples in future work. An application of these BV structures can be
found in \cite{M}. The notion of the BRST cohomology of the so-called
$W$-algebras has recently been investigated by Bouwknegt, McCarthy and Pilch
who again find a BV algebra structure (see \cite{BMP} and references therein).
The $W$-algebras can be viewed
as a generalization of the Virasoro algebra.

We note here that the abstract notions of a BV operator and BV algebra were
discovered quite independently of string theory by the mathematician J. L.
Koszul \cite{Koszul} (see also \cite{ks}). At roughly the same time Batalin and
Vilkovisky applied a particular BV operator to their ``antifield'' approach to
Lagrangian quantum field theory.  In a future work, we plan to review the
relationships between our own work and the relevant papers of Koszul, Batalin
and Vilkovisky, E. Getzler, A. Schwarz, and others.  We hope to highlight the
many intriguing parallel developments in physics and mathematics.  Such
parallels seem to be the hallmark of string theory itself.

{\bf Acknowledgments}: We thank G. Moore for many informative discussions about
the operator product expansions. We thank Y. Kosmann-Schwarzbach for sending us
her recent report, and F. Akman for carefully proofreading our manuscript.

\section{Quantum Operator Algebras}\lb{sec2}

Let $V$ be a {\it bounded} graded vector space. By this, we mean $V$ is a \bZ
doubly graded vector space $V=\oplus V^n[m]$ such that for each fixed $n$,
$V^n[m]=0$ for all but finitely many negative $m$'s.
The degrees of a homogeneous element $v$ in $V^n[m]$ will
be denoted by $|v|=n$, $||v||=m$ respectively.
In physical applications, $|v|$ will be the fermion or ghost number of $v$.
In conformal field theory, $||v||$ will be the conformal dimension or weight of
$v$.

Let $z$ be
a formal variable with degrees $|z|=0$, $||z||=-1$. Then
it makes sense to speak of a {\it homogeneous} (biinfinite)
formal power series
\eq
u(z)=\sum_{n\in\bf Z}u(n)z^{-n-1}
\en
of degrees $|u(z)|$, $||u(z)||$ where
the coefficients $u(n)$ are homogeneous linear maps in $V$ with degrees
$|u(n)|=|u(z)|$, $||u(n)||=-n-1+||u(z)||$. Note then that
the terms $u(n)z^{-n-1}$ indeed have the same degrees
$|u(z)|$, $||u(z)||$ for all $n$. We call a finite sum of such
homogeneous series $u(z)$ a {\it quantum operator} on $V$,
and we denote the linear space of quantum operators as $QO(V)$.

{\it Notations: By the expression $(z-w)^n$, $n$ an integer, we usually mean
its formal power series expansion in the region $|z|>|w|$. Thus $(z-w)^{-2}$
and $(-w+z)^{-2}$ are different, as power series. When such expressions
are to be regarded as rational functions rather than formal series, we will
explicitly mention so. When $A(z)=\sum A(n)z^{-n-1}$ is a formal series with
coefficients $A(n)$ in whatever linear space, we define $Res_zA(z)=A(0)$,
$A(z)^+=\sum_{n\geq0} A(n)z^{-n-1}$, $A(z)^-=\sum_{n<0} A(n)z^{-n-1}$,
$\partial A(z)=\sum -(n+1)A(n)z^{-n-2}$. If $u(z),u'(z)$ belong to QOAs $O,O'$
respectively, we abbreviate $u(z)\otimes u'(z)$, as an element of $O\otimes
O'$, simply as $u(z)u'(z)$. When no ambiguity occurs, we denote
$|u(z)|,||u(z)||$ simply as $|u|,||u||$. The restricted dual of a graded vector
space $V$ is denoted $V^\#$. If
$A_1(z)$,$A_2(z)$,... are quantum operators,
 an arbitrary matrix element
$\lgl x, A_1(z_1)A_2(z_2)\cdots y\rgl$ with $x\in V^\#$, $y\in V$, is denoted
as $\lgl A_1(z_1)A_2(z_2)\cdots\rgl$. In the interest of clarity, we often
write signs like $(-1)^{|t||u|}$ simply as $\pm$. This convention is used only
when the sign arises from permutation of elements. When in doubt, the reader
can easily recover the correct sign from such a permutation.}

Given two quantum operators $u(z),v(z)$, we write
\eq
:u(z)v(w):=u(z)^-v(w)+(-1)^{|u||v|}v(w)u(z)^+.
\en
Because $V$ is bounded, it's easy to check that if we replace $w$ by $z$, the
right hand side makes sense as a quantum operator and hence defines a
nonassociative product on $QO(V)$. It is called the {\it Wick product}.
Similarly given $u_1(z), \cdots,u_n(z)$, we define $:u_1(z_1)\cdots u_n(z_n):$
inductively as\\ $:u_1(z_1)(:u_2(z_2)\cdots u_n(z_n):):$.

\be{dfn}
For each integer $n$ we define
a product on $QO(V)$:
\eq
u(w)\circ_nv(w)=Res_z u(z)v(w)(z-w)^n
-(-1)^{|u||v|}Res_zv(w)u(z)(-w+z)^n.
\en
\end{dfn}

Explicitly we have:
\eq\lb{deriv}
u(z)\circ_n v(z)=\left\{\be{array}{ll}
\frac{1}{(-n-1)!}:\partial^{-n-1}u(z) \ v(z):& \mbox{if $n<0$}\\
{[(\sum_{m=0}^n
\left(\be{array}{c}
n\\
m
\end{array}\right)
u(m)(-z)^{n-m}),v(z)]} & \mbox{if $n\geq0$.}\end{array}\right.
\en
If $A$ is a  homogeneous linear operator on $V$, then it's clear that the
graded commutator $[A,-]$ is a graded derivation of each of the products
$\circ_n$.
Since $u(z)\circ_0 v(z) = [u(0),v(z)]$, we have
\be{pro}\lb{circle0}
For any $t(z),u(z),v(z)$ in $QO(V)$ and $n$ integer, we have
\[
t(z)\circ_0(u(z)\circ_n v(z)) = {[ t(z)\circ_0 u(z)]} \circ_n v(z)  \pm
u(z)\circ_n {[t(z)\circ_0 v(z)]},
\]
ie. $t(z)\circ_0$ is a derivation of every product in $QO(V)$.
\end{pro}

\be{pro}
For $u(z),v(z)$ in $QO(V)$, the following equality of
formal power series in two variables holds:
\eq\lb{ope}
u(z)v(w)=\sum_{n\geq0}u(w)\circ_n v(w) (z-w)^{-n-1}
+:u(z)v(w):.
\en
\end{pro}
Proof: We have $u(z)v(w)=[u(z)^+,v(w)]+:u(z)v(w):$. On the other hand by
inverting the second eqn. in \erf{deriv}, we  get
\eq
[u(m),v(w)]=\sum_{n=0}^m
\left(\be{array}{c}
m\\
n
\end{array}\right)
u(w)\circ_nv(w) w^{m-n}.
\en
Thus we have
\eqa
{[u(z)^+,v(w)]}&=&\sum_{m\geq n\geq0}
\left(\be{array}{c}
m\\
n
\end{array}\right)
u(w)\circ_nv(w) w^{m-n} z^{-m-1}\nnb\\
&=&\sum_{n\geq0}u(w)\circ_nv(w)\frac{1}{n!}\partial_w^n(z-w)^{-1}\nnb\\
&=&\sum_{n\geq0}u(w)\circ_n v(w) (z-w)^{-n-1}.\ \ \ \Box
\ena

In the sense of the above Proposition, $:u(z)v(w):$ is the {\it nonsingular
part} of
the {\it operator product expansion} \erf{ope}, while
$u(w)\circ_n v(w) (z-w)^{-n-1}$ is the {\it polar part}
of order $-n-1$ (see \cite{BPZ}).
The products $\circ_n$ will become important for describing
the algebraic and analytic structures of certain algebras
of quantum operators. Thus we introduce the following mathematical definitions:

\be{dfn}\lb{2.5}
A graded subspace \cA \ of $QO(V)$ containing the identity operator
and closed with respect to
all the products $\circ_n$ is called a quantum operator algebra.
 We say that $u(z)$ is local to $v(z)$
if $u(z)\circ_nv(z)=0$ for all but finitely many
positive $n$.
A QOA \cA is called local if its elements are pairwise mutually
local.
\end{dfn}
We observe that for any element $a(z)$ of a QOA, we have
$a(z)\circ_{-2}1 = \partial a(z)$. Thus a QOA is closed with respect
to formal differentiation.
\be{pro}\lb{2.6}
Let $u(z),v(z)$ be quantum operators, and $N$ a nonnegative integer. If
$u(z)\circ_nv(z)=0$ for $n\geq N$, then
$\lgl u(z)v(w)\rgl$ represents a rational function in $|z|>|w|$ with
poles along $z=w$ of order at most $N$.
\end{pro}
Proof: By eqn \erf{ope}, we have
\eq
\lgl u(z)v(w)\rgl = \sum_{n\geq0}\lgl u(w)\circ_n v(w)\rgl
(z-w)^{-n-1} +\lgl :u(z)v(w):\rgl.
\en
It is trivial to check that $\lgl :u(z)v(w):\rgl$, $\lgl u(w)\circ_n v(w)\rgl
\in \bC{[z^{\pm1},w^{\pm1}]}$. Thus our claim follows immediately.
$\Box$

\be{lem}\lb{operational}
Let $u(z)$ be local to $v(z)$, and $\lgl u(z)v(w)\rgl$ represent the
rational function $f(z,w)$. Then for $|w|>|z-w|$,
\eq
f(z,w)=\sum_{n\in{\bf Z}} \lgl u(w)\circ_nv(w)\rgl (z-w)^{-n-1}.
\en
\end{lem}
Proof: The Laurent polynomial $\lgl :u(z)v(w):\rgl$ in the above region is
just\\
$\sum_{i\geq0} \frac{1}{i!}\lgl :(\partial^i u(w))v(w):\rgl (z-w)^i$. Now apply
eqn. \erf{deriv}. $\Box$

We note that none of the products $\circ_n$ is associative
in general. However it clearly makes sense to speak of
the left, right or two sided ideals in a QOA as well as homomorphisms of QOAs
and they are defined in an obvious way. For example, a linear map $f:O\ra O'$
is a homomorphism if $f(u(z)\circ_n v(z))=fu(z)\circ_n fv(z)$ for all
$u(z),v(z)\in O$, and $f(1)=1$.

\be{dfn}\lb{2.6b}
Two quantum operators $u(z),v(z)$ are said to commute if they are mutually
local, and
$\lgl u(z)v(w)\rgl$,$\pm \lgl v(w)u(z)\rgl$
represent the same
rational function. This is equivalent (Proposition \ref{2.6}) to the following:
for some $N\geq0$,
$(z-w)^N\lgl u(z)v(w)\rgl= \pm(z-w)^N\lgl v(w)u(z)\rgl$
as Laurent polynomials. We call a QOA $O$ whose elements pairwise commute a
commutative QOA.
\end{dfn}

\be{pro}
If $u(z),v(z)$ commute, then for all $m$
\eq
[u(m),v(w)]=\sum_{n\geq0}
\left(\be{array}{c}
m\\
n
\end{array}\right)
u(w)\circ_nv(w) w^{m-n}.
\en
\end{pro}
Proof:
The case $m\geq0$ is obtained by inverting the second eqn. in \erf{deriv}.
Since $u(z)v(w)=[u(z)^+,v(w)]+:u(z)v(w):$ and
$v(w)u(z)=\mp[u(z)^-,v(w)]\pm:u(z)v(w):$, it follows from commutativity that
$\lgl [u(z)^-,v(w)]\rgl$ represents the same rational function as
$-\lgl [u(z)^+,v(w)]\rgl$ does, which is just
$-\sum_{n\geq0}\frac{u(w)\circ_nv(w)}{(z-w)^{n+1}}$. This gives
\eq
[u(z)^-,v(w)]=-\sum_{n\geq0}u(w)\circ_nv(w)(-w+z)^{-n-1}.
\en
Taking $Res_z [u(z)^-,v(w)]z^m$ for $m<0$ gives the desired result. $\Box$

The notion of commutativity here is closely related to the physicists' notion
of duality in conformal field theory\cite{MS}. Frenkel-Lepowsky-Meurman have
reformulated
the axioms of a VOA in terms of what they call rationality,
associativity and commutativity. The notion of commutativity in Definition
\ref{2.6b} is essentially the same as FLM's. This notion has also been
reformulated in the language of formal variables in \cite{dl}.

\be{dfn}
Let $\cO$ be a bounded graded space equipped with a distinguished vector $1$
and a set of bilinear products $\circ_n$, with $|1|=||1||=0$, $|\circ_n|=0,
||\circ_n||=-n-1$. We call $\cO$ a semi-inifinite commutative algebra if the
following holds: for homogeneous $u,v\in \cO$,\\
(i) $u\circ_n 1=\delta_{n,-1}u$ for $n\geq-1$;\\
(ii) $u\circ_n v=(-1)^{|u||v|}\sum_{p\in\bf
Z}(-1)^{p+1}(v\circ_pu)\circ_{n-p-1}1$.
\end{dfn}
Note that the sum in (ii) is finite because by (i), the summand is zero for
$p<n$, and by boundedness of the
space $\cO$, $v\circ_p u=0$ for $p>>0$. Note that the leading term on the RHS
of (ii) is $\pm v\circ_nu$. Thus the products $\circ_n$ are graded commutative
up to ``higher order corrections''. We claim that in a
$\frac{\infty}{2}$-commutative algebra, $1\circ_n t=\delta_{n,-1}t$ for all
$n$. Applying (i) to (ii) (with $u=1,v=t$), we have for $n\geq-1$,
\eq\lb{tempo}
1\circ_nt=\sum
(-1)^{p+1}(t\circ_{p}1)\circ_{n-p-1}1=(t\circ_{-1}1)\circ_n1=\delta_{n,-1}t.
\en
Now let $u=t,v=1$. Applying \erf{tempo}, (ii) becomes for $n\leq-2$,
\eqa
t\circ_n1&=&\sum(-1)^{p+1}(1\circ_pt)\circ_{n-p-1}1\nnb\\
&=&t\circ_{n}1+\sum_{p\leq-2}(-1)^{p+1}(1\circ_pt)\circ_{n-p-1}1.
\ena
Thus $\sum_{p\leq-2}(-1)^{p+1}(1\circ_pt)\circ_{n-p-1}1=0$. For $n=-2$, this
reads $0=-(1\circ_{-2}t)\circ_{-1}1=-1\circ_{-2}t$.
By induction on $n$, we have $1\circ_{n}t=0$ for $n\leq-2$.

\be{pro}
Let $O$ be a local QOA. Then $O$ is commutative iff it is a
$\frac{\infty}{2}$-commutative algebra.
\end{pro}
Proof: The unital property (i) follows from eqn. \erf{deriv}:
$u(z)\circ_n1=\delta_{n,-1}u(z)$ for $n\geq-1$.  Let $f,g$ be the rational
functions represented by $\lgl u(z)v(w)\rgl$ and $\lgl v(w)u(z)\rgl$
respectively. By Lemma \ref{operational}, we have
\eqa
f(z,w)&=&\sum_{n\in{\bf Z}} \lgl u(w)\circ_nv(w)\rgl (z-w)^{-n-1}\ \ \ \ for\
|w|>|z-w|\nnb\\
g(z,w)&=&\sum_{p\in{\bf Z}}\lgl v(z)\circ_pu(z)\rgl (w-z)^{-p-1}\ \ \ \ for\
|z|>|z-w|\nnb\\
&=&\sum_{p\in{\bf Z},m\geq0} (-1)^{p+1}\frac{1}{m!}
\lgl\partial^m (v(w)\circ_pu(w))\rgl (z-w)^{m-p-1}\ \ \ \ for\
|z|,|w|>|z-w|\nnb\\
&=&\sum_{p,m\in{\bf Z}} (-1)^{p+1}
\lgl (v(w)\circ_pu(w))\circ_{-m-1}1\rgl (z-w)^{m-p-1}\nnb\\
&=&\sum_{p,n\in{\bf Z}} (-1)^{p+1}
\lgl (v(w)\circ_pu(w))\circ_{n-p-1}1\rgl (z-w)^{-n-1}.
\ena
It's now clear that $f(z,w)=\pm g(z,w)$ holds iff the identity (ii) holds.
$\Box$

\section{Wick's calculus}

 In this section, we derive a number of
useful formulas relating various iterated products among three quantum
operators. Most of
these formulas are well-known to physicists who are familiar with the calculus
of operator product
expansions. We will also include a lemma on commutativity.

Let $t(z),u(z),v(z)$ be homogeneous quantum operators which pairwise commute.
\be{lem}(see \cite{li})\lb{lilemma}
For all $n$, $t(z)\circ_n u(z)$ and $v(z)$ commute.
\end{lem}
Proof: We include Li's proof here for completeness. For a positive integer $N$,
$(z-w)^{2N}$ is a binomial sum of terms $(z-x)^i(x-w)^{2N-i}$, $i=1,..,2N$.
So $(z-w)^{N+2N}(t(z)\circ_n u(z)) v(w)$ is a binomial sum of terms
\eq
Res_x\left((z-w)^N(z-x)^i(x-w)^{2N-i} (t(x)u(z)(x-z)^n\mp u(z)t(x)(-z+x)^n)v(w)
\right).
\en
We want to show that for large enough $N$, and for $0\leq i\leq 2N$, term by
term we have
\eqa\lb{comm}
&&(z-w)^N(z-x)^i(x-w)^{2N-i} (t(x)u(z)(x-z)^n\mp u(z)t(x)(-z+x)^n)v(w)\nnb\\
&&=\pm(z-w)^N(z-x)^i(x-w)^{2N-i} v(w)(t(x)u(z)(x-z)^n\mp u(z)t(x)(-z+x)^n).
\ena
Consider two cases: $i\geq N$ and $i<N$.
By assumption, $(z-x)^k(t(x)u(z)(x-z)^n\mp u(z)t(x)(-z+x)^n)=0$ for all large
enough $k$. So for large enough $N$, \erf{comm} holds for $i\geq N$.
Similarly for $i<N$, $(z-w)^N(x-w)^{2N-i} (t(x)u(z)(x-z)^n\mp
u(z)t(x)(-z+x)^n)v(w)$ coincides with  $\pm(z-w)^N(x-w)^{2N-i}
v(w)(t(x)u(z)(x-z)^n\mp u(z)t(x)(-z+x)^n)$. This shows that \erf{comm} holds
for each $i$. $\Box$

This lemma is useful for showing existence of commutative QOAs: it says that
given a set of pairwise commuting quantum operators, the QOA generated by the
set  is commutative. We now develop some abstract tools for studying the
structure of commutative QOAs.

Applying \erf{ope}, we have
\eqa\lb{tuvope}
&&:t(z)u(z):v(w) \nnb\\
&=&(t(z)^-u(z) \pm u(z)t(z)^+)v(w)\nnb\\
&=&t(z)^-u(z)v(w)\pm u(z)v(w)t(z)^+ \pm u(z)[t(z)^+,v(w)]\nnb\\
&=&\sum_{n\geq0}:t(z) (u(w)\circ_nv(w)):(z-w)^{-n-1}+
:t(z)u(z)v(w):\nnb\\
&&\pm\sum_{n,m\geq0}u(w)\circ_m(t(w)\circ_nv(w))(z-w)^{-n-m-2}+\nnb\\
&&\pm\sum_{n\geq0}:u(z) (t(w)\circ_nv(w)):(z-w)^{-n-1}.
\ena
Similarly,
\eqa\lb{tuvope2}
&&t(z):u(w)v(w):\nnb\\
&=&\pm[u(w)^-,t(z)]v(w)\pm u(w)^-t(z)v(w)\pm t(z)v(w)u(w)^+\nnb\\
&=&\pm\sum_{n,m\geq0}(-1)^{n+1}
 u(w)\circ_n(t(w)\circ_mv(w))(z-w)^{-n-m-2}\nnb\\
&&\pm\sum_{n\geq0}(-1)^{n+1}:u(z)\circ_n t(z)\ v(w):(z-w)^{-n-1}\nnb\\
&&\pm\sum_{n\geq0}:u(w)\ t(w)\circ_nv(w):(z-w)^{-n-1}
\pm:u(w)t(z) v(w):
\ena

\be{lem}\lb{2.6a}
The following
equalities hold in $|w|>|z-w|$:
\eqa
(i)&& \sum_{k\in {\bf
Z}}\frac{\lgl(:t(w)u(w):)\circ_kv(w)\rgl}{(z-w)^{k+1}}\nnb\\
&&=\sum_{n,m\geq0}\frac{\lgl:\partial^mt(w)\ u(w)\circ_nv(w):\rgl
\pm\lgl:\partial^m u(w)\ t(w)\circ_nv(w):\rgl}{m!(z-w)^{n-m+1}}\nnb\\
&&\pm\sum_{n,m\geq0}
\frac{\lgl u(w)\circ_n(t(w)\circ_mv(w))\rgl}{(z-w)^{n+m+2}}\nnb\\
&&+\sum_{m\geq0}\frac{\lgl:\partial^m (t(w)u(w))\ v(w):\rgl}{m!(z-w)^{-m}}\\
(ii)&&\pm\sum_{k\in {\bf Z}}\frac{\lgl
t(w)\circ_k:u(w)v(w):\rgl}{(z-w)^{k+1}}\nnb\\
&&=\sum_{n,m\geq0}(-1)^{n+1}
\frac{\lgl u(w)\circ_n(t(w)\circ_mv(w))\rgl}{(z-w)^{n+m+2}}\nnb\\
&&+\sum_{n,m\geq0}(-1)^{n+1}\frac{\lgl:\partial^m(u(w)\circ_n t(w))v(w):\rgl}
{m!(z-w)^{n-m+1}}\nnb\\
&&+\sum_{n\geq0}\frac{\lgl:u(w)\ t(w)\circ_nv(w):\rgl}{(z-w)^{n+1}}\nnb\\
&&+\sum_{m\geq0}\frac{\lgl :u(w)(\partial^m t(w)) v(w):\rgl}{m!(z-w)^{-m}}
\ena
\end{lem}
Proof:
 To prove (i), consider  matrix coefficients on both sides of  eqn.
\erf{tuvope}. By assumption of
commutativity these matrix coefficients represent  rational functions.
Expanding both sides using Lemma \ref{operational}, we get the first eqn. (i).
The eqn. (ii) is derived similarly from \erf{tuvope2}. $\Box$

By reading off coefficients of the $(z-w)^i$, we can use this lemma to
simultaneously compute all products $:t(w)u(w):\circ_k v(w)$, and
$t(w)\circ_k :u(w)v(w):$ in terms of the products among the constituents
$t(w),u(w),v(w)$. Thus it is a kind of recursion relation for the products.
In the examples below, we will see how it allows us to understand the structure
of commutative QOAs.

\be{lem}\lb{abelianQOA}
If $t(z)^\pm u(w)^\pm=(-1)^{|t||u|}u(w)^\pm t(z)^\pm$, then
$:t(z)u(w)v(x):=(-1)^{|t||u|}:u(w)t(z)v(x):$.
\end{lem}
Proof: Applying the definition of the Wick product (and surpressing $z,w,x$):
\eqa
&&:tuv:-(-1)^{|t||u|}:utv:\nnb\\
&&=t^-(u^-v+(-1)^{|u||v|}vu^+)+
(-1)^{|t|(|u|+|v|)}(u^-v+(-1)^{|u||v|}vu^+)t^+\nnb\\
&&-(-1)^{|t||u|}\left(u^-(t^-v+(-1)^{|t||v|}vt^+)+
(-1)^{|u|(|y|+|v|)}(t^-v+(-1)^{|t||v|}vt^+)u^+\right)\nnb\\
&&=0.\ \ \ \Box
\ena

\subsection{Examples}

Let $QO(V)^-=\{u(z)^-\ |\ u(z)\in QO(V)\}$. This space is obviously closed
under differentiation and the Wick product. It follows that the space is also
closed under all $\circ_n$, $n$ negative. Also observe that for any
$u(z),v(z)\in QO(V)$, we have $u(z)^-v(w)^-=:u(z)^- v(w)^-:$. It follows that
the products $\circ_n$, $n=0,1,...$, restricted to $QO(V)^-$, all vanish. Thus
$QO(V)^-$ is a local QOA.

Let $LO(V)$ be the algebra of graded linear operators on $V$. We can regard
each operator $A$ as a formal series with just the constant term. This makes
$LO(V)$ a subspace of $QO(V)$. It is obvious that every $\circ_n$ restricted to
$LO(V)$ vanishes except for $n=-1$, in which case $\circ_{-1}$ is the usual
product on $LO(V)$. Thus $LO(V)$ is a very degenerate example of a QOA.
Obviously, any commutative subalgebra of $LO(V)$ is a commutative QOA.

Let $\cC$ be the Clifford algebra with the generators $b(n),c(n)$
($n\in\bZ$) and the relations \cite{FMS}\cite{KR}\cite{Akman}
\be{eqnarray}
b(n)c(m)+c(m)b(n)&=&\delta_{n,-m-1}\nnb\\
b(n)b(m)+b(m)b(n)&=&0\nnb\\
c(n)c(m)+c(m)c(n)&=&0
\end{eqnarray}
Let $\lambda$ be a fixed integer.
The algebra \cC becomes \bZ-bigraded if we define the degrees
$|b(n)|=-|c(n)|=-1$, $||b(n)||=\lambda-n-1$, $||c(n)||=-\lambda-n$.
Let $\bigwedge^*$ be the graded irreducible
$\cC^*$-module with generator \bone\ and relations
\be{equation}
b(m)\bone=c(m)\bone=0, \ \ \ m\geq0
\end{equation}

Let $b(z),c(z)$ be the quantum operators
\be{eqnarray}
b(z)&=&\sum_mb(m)z^{-m-1}\nnb\\
c(z)&=&\sum_mc(m)z^{-m-1}
\end{eqnarray}
Let $O(b,c)$ be the smallest QOA containing $b(z),c(z)$.

\be{pro}\lb{Obc}
The QOA $O(b,c)$ is commutative. It has a basis consisting of the monomials
\eq\lb{bcmono}
:\partial^{n_1}b(z)\cdots\partial^{n_i}b(z)\
\partial^{m_1}c(z)\cdots\partial^{m_j}c(z):
\en
with $n_1>...>n_i\geq0$, $m_1>...>m_j\geq0$.
\end{pro}
Proof: Computing the OPE of $b(z),c(w)$, we have
\eqa\lb{bcope}
b(z)c(w)&=&(z-w)^{-1}+:b(z)c(w):\nnb\\
c(w)b(z)&=&(w-z)^{-1}+:c(w)b(z):\nnb\\
:b(z)c(w):&=&-:c(w)b(z):.
\ena
It follows that $b(z)$ and $c(z)$ commute. Also $b(z),c(z)$ each commutes with
itself, hence they form a pairwise commuting set.
By Lemma \ref{lilemma}, they generate a commutative QOA.

If each $u_1(z),...,u_k(z)$ is
of the form  $\partial^{n}b(z)$ or $\partial^{m}c(z)$, let's
call $:u_{1}(z)\cdots u_{k}(z):$ a monomial of degree $k$.
We claim that it's proportional to some monomial \erf{bcmono} with
$n_1>...>n_i\geq0$, $m_1>...>m_j\geq0$. If $t(z), u(z)$ each is of the form
$\partial^{n}b(z)$ or $ \partial^{m}c(z)$, it is easy to check that
$t(z)^\pm u(z)^\pm=-u(z)^\pm t(z)^\pm$. It follows from Lemma \ref{abelianQOA}
that
$:t(z)u(z)v(z):=-u(z)t(z)v(z):$ for any element $v(z)\in O(b,c)$. This shows
that $:u_{1}(z)\cdots u_{k}(z):$ is equal to
$(-1)^\sigma:u_{\sigma(1)}(z)\cdots u_{\sigma(k)}(z):$ for any permutation
$\sigma$ of $1,...,k$.

Let $O'$ be the linear span of the monomials \erf{bcmono}.
We now show that $A\circ_k B\in O'$ for any $k$ and any two monomials $A,B$,
hence $O(b,c)=O'$.
We will do a double induction on the degrees of $A$ and $B$.
Case 1: let $A=t(z)$, $B=:u(z)v(z):$ with
$t(z), u(z)$ each monomial of degree 1, and
$v(z)$ of any degree.
If $v(z)=1$, then by \erf{bcope} $t(w)\circ_k:u(w)v(z):\in O'$.
 By induction on the degree of $v(z)$ and applying
Lemma \ref{2.6a}(ii), we see that $t(w)\circ_k:u(w)v(w):\in O'$. This shows
that $A\circ_kB\in O'$ for $A$ of degree 1, $B$ of any degree.
Now case 2: suppose $A=:t(z)u(z):$, $B=v(z)$, where $t(z)$ is of degree 1 and
$u(z),v(z)$ of any degree. By induction on the degree of $u(z)$, it's clear
from Lemma \ref{2.6a}(i) that this case reduces to case 1.

Finally we must show that the monomials \erf{bcmono} are linearly independent.
 Define a map $O(b,c)\ra \bigwedge$ by $u(z)\mapsto u(-1)\bone$. We see that
this map gives a 1-1 correspondence between the set of monomials \erf{bcmono}
and a basis of $\bigwedge$. This completes the proof. $\Box$.

Let $M(\kappa,0)$ be the Verma module of the Virasoro algebra with highest
weight $(\kappa,0)$ and vacuum vector $v_0$. Let $M(\kappa)$ be the quotient of
$M(\kappa,0)$ by the submodule generated by $L_{-1}v_0$. Let $O_\kappa(L)$ be
the QOA generated by $L(z)=\sum L_nz^{-n-2}$ in $QO(M(\kappa))$.
\be{pro}\lb{Okappa}(see \cite{BPZ}\cite{FZ})
The QOA $O_\kappa(L)$ is commutative. It has a basis consisting of monomials
\eq\lb{Lmono}
:\partial^{n_1}L(z)\cdots\partial^{n_i}L(z):
\en
with $n_1\geq...\geq n_i\geq0$.
\end{pro}
Proof: A direct computation gives
\eqa\lb{LOPE}
{[L(z)^+,L(w)]}&=&\frac{\kappa}{2}(z-w)^{-4}+2L(w)(z-w)^{-2}+\partial
L(w)(z-w)^{-1}\nnb\\
{[L(z)^-,L(w)]}&=&-\frac{\kappa}{2}(w-z)^{-4}-2L(w)(w-z)^{-2}+\partial
L(w)(w-z)^{-1}.
\ena
But we also have $L(z)L(w)=[L(z)^+,L(w)]+:L(z)L(w):$, and
$L(w)L(z)=-[L(z)^-,L(w)]+:L(z)L(w):$. Combining these with
\erf{LOPE}, it is obvious that $\lgl L(z)L(w)\rgl$ and $\lgl L(w)L(z)\rgl$
represent the same rational function.
Thus $L(z)$ commutes with itself as a quantum operator. By Lemma \ref{lilemma},
$O_\kappa(L)$ is commutative.

Let $O'$ be the linear span of the monomials \erf{Lmono} with $n_1,..,n_i\geq0$
{\it unrestricted}. To show $O'$ is closed under all the products (hence
$O_\kappa(L)=O'$), we apply induction and Lemma
\ref{2.6a} as in the case of $O(b,c)$ above. We now show that
we can restrict to those monomials \erf{Lmono} with $n_1\geq,...\geq n_i\geq0$,
and that the resulting monomials form
a basis. First by direct computation, we see
that $O_\kappa(L)$ is a $Vir$-module defined by the action
\eq
L(n)\cdot u(z)=L(z)\circ_n u(z).
\en
Since $O_\kappa(L)$ is spanned by the monomials \erf{Lmono}, and because
$L(-n-1)\cdot u(z)=\frac{1}{n!}:\partial^{n} L(z)u(z):$ for $n\geq0$, it
follows that the module is
cyclic. Thus we have a unique onto map of $Vir$-modules $M(\kappa)\ra
O_\kappa(L)$ sending $v_0$
to $1$. But $M(\kappa)$ has a PBW basis consisting of $L(-n_1-1)\cdots
L(-n_i-1)v_0$,
$n_1\geq,...\geq n_i\geq0$. This shows that the monomials \erf{Lmono} with
$n_1\geq,...\geq n_i\geq0$
span $O_\kappa(L)$. Now define a map $O_\kappa(L)\ra M(\kappa)$ by $u(z)\mapsto
u(-1)v_0$. This is the inverse to the previous map, hence it must map a basis
to a basis. $\Box$

\section{BRST cohomology algebras}

\be{dfn}
A conformal QOA with central charge $\kappa$ is a pair $(O,f)$, where $O$ is a
commutative QOA equipped with a homomorphism $f:O_\kappa(L)\ra O$ such that for
every homogeneous $u(z)\in O$,
\eq
fL(z)u(w) =\cdots + ||u||u(w)(z-w)^{-2} + \partial u(w) (z-w)^{-1}+:fL(z)u(w):
\en
where ``$\cdots$'' denotes the higher order polar terms.  In other words,
$fL(z)\circ_1 u(z)=||u||u(w)$ and $fL(z)\circ_0 u(z)=\partial u(z)$.  For
simplicity, we sometimes write $f:O_\kappa(L)\ra O$, or simply $O_\kappa(L)\ra
O$, to denote a conformal
QOA. A homomorphism $(O,f)\stackrel{h}{\ra}(O',f')$ of conformal QOAs is a
homomorphism of QOAs $h:O\ra O'$ such that $h\circ f=f'$.
\end{dfn}

\be{lem}
Let $O$ be a commutative QOA generated by a set $S$. Let $X(z)\in O$ such
that for all $u(z)\in S$,
\eqa\lb{Xvir}
X(z)u(w) &=&\cdots + ||u||u(w)(z-w)^{-2} + \partial u(w)
(z-w)^{-1}+:X(z)u(w):\nnb\\
X(z)X(w) &=&\frac{\kappa}{2}(z-w)^{-4}+2X(w)(z-w)^{-2} + \partial X(w)
(z-w)^{-1}+:X(z)X(w):.
\ena
Then there's a unique homomorphism $f:O_\kappa(L)\ra O$ such that $fL(z)=X(z)$.
Moreover $(O,f)$ is a conformal QOA with central charge $\kappa$.
\end{lem}
Proof: Uniqueness of $f$ is clear because $O_\kappa(L)$ is generated by $L(z)$.

Step 1: We claim that $O$ is a $Vir$-module defined by the action:
\eq
L(n)\cdot u(z)=X(z)\circ_n u(z).
\en
We must show that
\eq\lb{vircomm}
L(m)\cdot L(n)\cdot u(z)-L(n)\cdot L(m)\cdot u(z)
=(m-n)L(m+n-1)\cdot u(z)+\frac{\kappa}{12}m(m-1)(m-2)\delta_{n+m-2}u(z).
\en
{}From the OPE of $X(z)$, we get (cf. \erf{LOPE}):
\eqa
LHS&=&Res_{z_1}Res_{z_2}{[X(z_2),X(z_1)]}u(z)(z_2-z)^m(z_1-z)^n\nnb\\
&&-Res_{z_1}Res_{z_2}u(z){[X(z_2),X(z_1)]}(-z+z_2)^m(-z+z_1)^n\nnb\\
&=&Res_{z_1}Res_{z_2}(z_2-z)^m(z_1-z)^n\left(\frac{\kappa}{12}\partial_{z_1}^3+
2X(z_1)\partial_{z_1}+\partial_{z_1}X(z_1)\right)\delta(z_1,z_2)u(z)\nnb\\
&&-Res_{z_1}Res_{z_2}u(z)(-z+z_2)^m(-z+z_1)^n
\left(\frac{\kappa}{12}\partial_{z_1}^3+
2X(z_1)\partial_{z_1}+\partial_{z_1}X(z_1)\right)\delta(z_1,z_2)u(z)\nnb\\
&&
\ena
where $\delta(z_1,z_2)=(z_2-z_1)^{-1}-(-z_1+z_2)^{-1}$. For any Laurent
polynomial $g(z_1,z_2)$, and any formal series $h(z_1,z_2)$,
 we have the identities
\eqa
g(z_1,z_2)\delta(z_1,z_2)&=&g(z_1,z_1)\delta(z_1,z_2)\nnb\\
Res_{z_1}\delta(z_1,z_2)&=&1\nnb\\
Res_{z_1}g(z_1,z_2)\partial_{z_1}^k h(z_1,z_2)
&=&(-1)^kRes_{z_1}(\partial_{z_1}^kg(z_1,z_2))\ h(z_1,z_2).
\ena
Applying these and continuing the above computation:
\eqa
LHS&=&Res_{z_1}\left(-\frac{\kappa}{12}n(n-1)(n-2)(z_1-z)^{m+n-3}\right.\nnb\\
&&\left. -2X(z_1) n (z_1-z)^{m+n-1}-2\partial_{z_1}X(z_1)(z_1-z)^{m+n}
+\partial_{z_1}X(z_1)(z_1-z)^{m+n}\right) u(z)\nnb\\
&&-Res_{z_1}\left(-\frac{\kappa}{12}n(n-1)(n-2)(-z+z_1)^{m+n-3}\right.\nnb\\
&&\left. -2X(z_1) n (-z+z_1)^{m+n-1}-2\partial_{z_1}X(z_1)(-z+z_1)^{m+n}
+\partial_{z_1}X(z_1)(-z+z_1)^{m+n}\right) u(z)\nnb\\
&=&-\frac{\kappa}{12}n(n-1)(n-2)1\circ_{m+n-3}u(z)+(m-n)X(m+n-1)\cdot u(z).
\ena
This proves \erf{vircomm}.
Thus we've a map $f:O_\kappa(L)\cong M(\kappa)\ra O$ with $1\mapsto 1$.
Moreover,
we have
\eq\lb{fmono}
f(:\partial^{n_1}L(z)\cdots\partial^{n_i}L(z):)=
:\partial^{n_1}X(z)\cdots\partial^{n_i}X(z):.
\en

Step 2: We'll show that $f$ is a QOA homomorphism, ie.
 for all $A,B\in O_\kappa(L)$ and integers $n$,
\eq\lb{homo}
f(A\circ_nB)=fA\circ_n fB.
\en
It's enough to do it for $A,B$ being monomials $:u_1(z)\cdots u_k(z):$ where
each $u$ has degree 1, ie. of the form $\partial^nL(z)$. Once again by double
induction on the degrees of $A,B$, it's easy to show that $A\circ_n B$ can be
reduced, by {\it Lemma \ref{2.6a}} and the OPE $L(z)L(w)$ only, to a linear sum
of the above monomials. By \erf{fmono},  $fA\circ_n fB$ must be reduced, by
{\it Lemma \ref{2.6a}} and the OPE $X(z)X(w)$ only, to a linear
sum {\it identical} to the reduction of $f(A\circ_nB)$. This is so because
$L(z)$ and $X(z)$ have identical OPEs. This proves \erf{homo} for all $A,B$.

Step 3: We must show that $fL(z)=X(z)$ has the desired properties: $X(z)\circ_0
u(z)=\partial u(z)$ and $X(z)\circ_1 u(z)=||u|| u(z)$, for all $u(z)\in O$.
By assumption, they hold for all $u(z)$ in the generating set $S$. Checking the
properties is an easy exercise applying the formulas $X(z)\circ_0 \partial
u(z)=[X(0),u(z)]$, $X(z)\circ_1 u(z)=[X(1)-X(0)z,u(z)]$ and applying induction.
$\Box$

Consider, as an example, $O(b,c)$. For a fixed $\lambda$, let
\eq
X(z)=(1-\lambda):\partial b(z)\ c(z)-\lambda:b(z)\partial c(z):
\en
and $S=\{ b(z),c(z)\}$. Then we have, by direct computation \cite{FMS},
\eqa
X(z)b(w) &=&\lambda b(w)(z-w)^{-2} + \partial b(w) (z-w)^{-1}+:X(z)b(w):\nnb\\
X(z)c(w) &=&(1-\lambda) c(w)(z-w)^{-2} + \partial c(w)
(z-w)^{-1}+:X(z)c(w):\nnb\\
X(z)X(w) &=&\frac{\kappa}{2}(z-w)^{-4}+2X(w)(z-w)^{-2} + \partial X(w)
(z-w)^{-1}+:X(z)X(w):.
\ena
where $\kappa=-12\lambda^2+12\lambda-2$. It follows that we have a homomorphism
$f_\lambda: O_\kappa(L)\ra O(b,c)$ such that $(O(b,c),f_\lambda)$ is a
conformal QOA with central charge $\kappa$.

Similarly the pair $(O_\kappa(L),id)$ is itself a conformal QOA. Thus by
definition, it is the initial object in the category of conformal QOAs with
central charge $\kappa$.

\subsection{The BRST construction}

It is evident that if $(O,f)$, $(O',f')$ are conformal QOAs
on the respective spaces $V$, $V'$ with
central charges $\kappa,\kappa'$, then
$(O\otimes O',f\otimes f')$ is a conformal QOA on $V\otimes V'$
with central charge $\kappa+\kappa'$. From now on we fix $\lambda=2$ which
means that $(O(b,c),f_\lambda)$ now has central charge -26. Let $(O,f)$
be any conformal QOA with central charge $\kappa$ and consider
\eq
C^*(O)= O(b,c)\otimes O
\en
where $*$ denotes the total first degree.
For simplicity, we write $L^C(z)=f_\lambda L(z)+fL(z)$.
\be{pro}
For every conformal QOA $O$, there is a unique homogeneous element $J_O(z)\in
C^*(O)$ with the following properties:\\
(i) (Cartan identity) $J_O(z)b(w)=L^C(w)(z-w)^{-1}+:J_O(z)b(w):$.\\
(ii) (Universality) If $(O,f)\ra(O',f')$ is a homomorphism of conformal QOAs,
then the induced homomorphism $C^*(O)\ra C^*(O')$ sends $J_O(z)$ to
$J_{O'}(z)$.
\end{pro}
Proof: Since the category of conformal QOAs with central charge $\kappa$ has
$(O_\kappa(L),id)$ as the initial object, if we can show that there is a unique
$J_{O_\kappa(L)}$ satisfying property (i), then (ii) implies that the same
holds for every other object in that category.

Property (i) implies $|J_{O_\kappa(L)}|=1=||J_{O_\kappa(L)}||$. Let's list a
basis of $C^1(O_\kappa(L))[1]$ given by Propositions \rf{Obc}, \rf{Okappa}:
$:c(z)L(z):$, $:b(z)c(z)\partial c(z):$, $\partial^2c(z)$. Take a linear
combination of these elements and compute its OPE with $b(z)$. Now requiring
property (i), we determine the coefficients of the linear combination and get
\eq
J_{O_\kappa(L)}(z)=:c(z) L(z):+:b(z)c(z)\partial c(z):.
\en
Now given a conformal QOA $(O,f)$, the induced map $f^*:C^*(O_\kappa(L))\ra
C^*(O)$ sends $J_{O_\kappa(L)}(z)$ to $J_O(z)=:c(z) fL(z):+:b(z)c(z)\partial
c(z):$. This completes our proof. $\Box$.

It follows from property (i) that
\eq
L^C(z)=J_O(z)\circ_0b(z)=[Q,b(z)]
\en
where $Q=Res_z J_O(z)$.
\be{lem}\cite{KO}\cite{Fe}\cite{FGZ}
Let $(O,f)$ be a conformal QOA with central charge $\kappa$. Then $Q^2=0$ iff
$\kappa=26$.
\end{lem}
Proof: We'll drop the subscript for $J_O$ and write $fL(z)$ as $L(z)$. Let's
compute $2Q^2=[Q,Q]=Res_w[Q,J(w)]=Res_w J(w)\circ_0 J(w)$. Since $J(w)\circ_0
J(w)$ is the coefficient of $(z-w)^{-1}$ in the OPE $J(z)J(w)$, we can
extract this term from the OPE. Now $J(z)J(w)$ is
the sum of 4 terms:
\eqa
&(i)& c(z) L(z)c(w)L(w)\nnb\\
&(ii)& c(z)L(z):b(w)c(w)\partial c(w):\nnb\\
&(iii)& :b(z)c(z)\partial c(z): c(w)L(w)\nnb\\
&(iv)& :b(z)c(z)\partial c(z)::b(w)c(w)\partial c(w):.
\ena
 Extracting the coefficient of $(z-w)^{-1}$ (which is done by applying Lemma
\ref{2.6a} repeatedly) in each of these 4 OPEs, we
get respectively (surpressing $w$):
\eqa
&(i)& 2\partial c\ cL+\frac{\kappa}{12}\partial^3c\ c\nnb\\
&(ii)& c\partial c\ L\nnb\\
&(iii)& c\partial c\ L\nnb\\
&(iv)& \frac{3}{2}\partial(\partial^2 c\ c)-\frac{13}{6}\partial^3 c\ c
\ena
Thus $J(w)\circ_0J(w)=\frac{3}{2}\partial(\partial^2 c(w)\
c(w))+\frac{\kappa-26}{12}\partial^3 c(w)\ c(w)$.
The $Res_w$ of this is zero iff $\kappa=26$. $\Box$
%
%\eqa
%&(i)&c(z) L(z)c(w)L(w)=c(z)c(w)\left(\frac{\kappa}{12}(z-w)^{-4}+
%2L(w)(z-w)^{-2}+\partial L(w)(z-w)^{-1} +\cdots\right)\nnb\\
%&(ii)&c(z)L(z):b(w)c(w)\partial c(w):=L(z)c(w)\partial c(w)
%%%(z-w)^{-1}+\cdots\nnb\\
%&(iii)&:b(z)c(z)\partial c(z): c(w)L(w)=L(w)c(z)\partial
%%%c(z)(z-w)^{-1}+\cdots\nnb\\
%&(iv)&:b(z)c(z)\partial c(z)::b(w)c(w)\partial c(w):
%&&=b(z)^-c(z)\partial c(z):b(w)c(w)\partial c(w):
%-c(z)\partial c(z):b(w)c(w)\partial c(w):b(z)^+\nnb\\
%&&-c(z)\partial c(z):b(w)\partial c(w):(z-w)^{-1}
%+c(z)\partial c(z):b(w)c(w):(z-w)^{-2}\nnb\\
%&&=-:b(z)\partial c(z)c(w)\partial c(w):(z-w)^{-1}
%-:b(z)c(z)c(w)\partial c(w):(z-w)^{-2}
%:b(z)\partial c(z)c(w)\partial c(w):(z-w)^{-1}
%\ena
%

Recall that (Lemma \ref{circle0}) $[Q,-]=J(z)\circ_0$ is a derivation of the
QOA $C^*(O)$. For $\kappa=26$, which we assume from now on, $[Q,-]$ becomes
a differential on $C^*(O)$ and we have a cochain complex
\eq
[Q,-]:C^*(O)\lra C^{*+1}(O).
\en
It is called {\it the BRST complex} associated to $O$. Its cohomology will be
denoted as $H^*(O)$. All the operations $\circ_n$ on $C^*(O)$
descend to the cohomology. However, all but one is trivial.
\be{thm}\lb{5.3}\cite{Wi3}\cite{WZ}\cite{LZ9}
The Wick product $\circ_{-1}$ induces a graded commutative associative product
on
$H^*(O)$ with unit element represented by the identity operator. Moreover,
every cohomology class is represented by a quantum operator $u(z)$ with
$||u||=0$.
\end{thm}
Proof: Let $u(z),v(z)$ be two elements of $C^*(O)$ annihilated by $[Q,-]$.
Consider the rational functions represented by $\lgl u(z)v(w)\rgl$, $\pm\lgl
v(w)u(z)\rgl$. Expand the second one in the region $|w|>|z-w|$:
\eqa\lb{vu}
\lgl v(w)u(z)\rgl&=&\sum_{m\geq0}\lgl v(z)\circ_mu(z)\rgl (w-z)^{-m-1} +\lgl
:v(w)u(z):\rgl\nnb\\
&=&\sum_{m\geq0}\sum_{k\geq0} \partial^k(\lgl v(w)\circ_mu(w)\rgl)
\frac{(-1)^k}{k!} (w-z)^{k-m-1} \nnb\\
& &+\sum_{k\geq0} \lgl :v(w)\partial^ku(w):\rgl \frac{(-1)^k}{k!}(w-z)^k.
\ena
By commutativity, this coincides with the following rational function times
$(-1)^{|v||u|}$:
\eq
\lgl u(z)v(w)\rgl=\sum_{n\geq0}\lgl u(w)\circ_nv(w)\rgl (z-w)^{-n-1} +\lgl
:u(z)v(w):\rgl.
\en
Expanding this in the same region $|w|>|z-w|$, we see that its coefficient of
$(z-w)^0$ is $\lgl :u(w)v(w):\rgl$. Equating this with the same coefficient in
\erf{vu} (with appropriate sign), we get
\eq
\sum_{m\geq0}\frac{\partial^{m+1}(\lgl
v(w)\circ_mu(w)\rgl)}{(m+1)!}+\lgl:v(w)u(w):\rgl
=(-1)^{|v||u|}\lgl :u(w)v(w):\rgl.
\en
We can write the first term on the left hand side as $\lgl\partial A(w)\rgl$
where $A(w)$ is a quantum operator.
(We can do this because the sum is finite by locality.) Note that because
$[Q,-]$ is a derivation on the products $\circ_m$ and because $[Q,-]$
annihilates all $\partial^i u(w),\partial^jv(w)$, it follows that $[Q,A(w)]=0$.
Thus we have
\eqa
\lgl:v(w)u(w):\rgl-(-1)^{|v||u|}\lgl :u(w)v(w):\rgl.&=&\lgl\partial
A(w)\rgl\nnb\\
&=&\lgl L^C(w)\circ_0 A(w)\rgl\nnb\\
&=&\lgl[Q,b(w)]\circ_0 A(w)\rgl\nnb\\
&=&\lgl[Q,b(w)\circ_0A(w)]\rgl.
\ena
This implies that $:v(w)u(w):$, $(-1)^{|v||u|} :u(w)v(w):$ are cohomologous to
each other.

To prove associativity, let $u(z),v(z),t(z)$ represent three cohomology
classes. We will compute $\lgl :u(z)v(z): t(w)\rgl$ in two different ways.
First we have
\eqa
\lgl :u(z)v(z): t(w)\rgl&=&\lgl (u(z)^-v(z)+(-1)^{|u||v|} v(z)u(z)^+)
t(w)\rgl\nnb\\
&=&\sum_{n\geq0}\lgl :u(z)v(w)\circ_n t(w):\rgl(z-w)^{-n-1} +
\lgl:u(z)v(z)t(w):\rgl \nnb\\  & &+
(-1)^{|u||v|}\sum_{m\geq0}\lgl v(z)u(w)\circ_n t(w)\rgl(z-w)^{-m-1}.
\ena
The right hand side is now a rational function which we can expand in the
region $|w|>|z-w|$ and extract the coefficient of $(z-w)^0$. On the other hand
this coefficient must be equal to $\lgl :(:u(w)v(w):) t(w):\rgl$. Thus we get
\eqa
\lgl :(:u(w)v(w):)t(w):\rgl&=&\lgl :u(w)(:v(w)t(w):):\rgl\nnb\\
&+ &\sum_{n\geq0}\frac{\lgl:\partial^{n+1}u(w)\ v(w)\circ_nt(w):\rgl}{(n+1)!}
\nnb\\ & + &(-1)^{|u||v|}\sum_{m\geq0}\frac{\lgl :\partial^{m+1}v(w)\
u(w)\circ_n t(w):\rgl}{(m+1)!}.
\ena
Using the fact that $\partial A(w)=[Q,b(w)]\circ_0A(w)$ for any $A(w)$, and the
fact that $[Q,-]$ annihilates $u(w),v(w),t(w)$, we see that
$:(:u(w)v(w):)t(w):$ and $:u(w)(:v(w)t(w):):$ are cohomologous to each other.

Checking that the quantum operator $1$ represents the unit is a trivial
exercise. Finally we want to show that an element $u(z)$, with $||u||\neq0$,
annihilated by $[Q,-]$ is cohomologous to zero. Recall that
$L^C(z)\circ_1u(z)=||u|| u(z)$. But the left hand side is $[Q,b(z)\circ_1
u(z)]$ which is cohomologous to zero. $\Box$

\section{Batalin-Vilkovisky Algebras}

Let $A^*$ be a \bZ graded commutative associative algebra.
For every $a\in A$, let $l_a$ denote the linear map on $A$  given by
the left multiplication by $a$. Recall that a (graded) derivation $d$ on $A$ is
a homogeneous linear operator such that ${[d,l_a]} - l_{da}=0$ for all $a$.
A BV operator \cite{W2}\cite{Sch}\cite{GJ} $\Delta$ on $A^*$ is
a linear operator of degree -1 such that:\\
(i) $\Delta^2=0$;\\
(ii) ${[\Delta,l_a]} - l_{\Delta a}$ is a derivation on $A$ for all $a$, ie.
$\Delta$ is
a {\it second} order derivation.

A BV algebra is a pair $(A,\Delta)$ where $A$ is a graded commutative
algebra and $\Delta$ is a BV operator on $A$. The following is an elementary
but fundamental lemma:
\be{lem}\cite{Koszul}\cite{GJ}\cite{Pen}
Given a BV algebra $(A,\Delta)$, define the BV bracket $\{,\}$ on $A$ by:
\[
(-1)^{|a|}\{a,b\} = {[\Delta,l_a]}b - l_{\Delta a}b.
\]
Then $\{,\}$ is a graded Lie bracket on $A$ of degree -1.
\end{lem}
By property (ii) above, it follows immediately that for every
$a\in A$, $\{a,-\}$ is a derivation on $A$. Thus a BV algebra is a special case
of an odd Poisson algebra which, in mathematics, is also known as a
Gerstenhaber
algebra \cite{Gers1}\cite{Gers2}. We note that $A^1$ is canonically a Lie
algebra.

Consider now the linear operator $\Delta: C^*(O)\lra C^{*-1}(O)$,
$u(z)\mapsto b(z)\circ_1 u(z)$.
\be{thm}\lb{5.4}\cite{LZ9}
The operator $\Delta$ descends to the cohomology $H^*(O)$. Morever,
it is a BV operator on the commutative algebra $H^*(O)$. Thus
$H^*(O)$ is naturally a BV algebra.
\end{thm}
Proof: By the theorem above, it is enough consider the action of $\Delta$ on
elements $u(z)$ with $||u||=0$. If $[Q,u(z)]=0$, we have
\eq
[Q,\Delta u(z)]=[Q,b(z)]\circ_1 u(z)=L^C(z)\circ_1u(z)=||u(z)||u(z)=0.
\en
Thus $\Delta$ is well-defined on the cohomology.

Also we have
\eq
\Delta^2u(z)=[b(1)-b(0)z,[b(1)-b(0)z,u(z)]]=0
\en
because the $b's$ anticommute.

We define the following bilinear operation:
\eq
(-1)^{|u|}\{u(z),v(z)\}=\Delta(:u(z)v(z):)-:(\Delta u(z))v(z):-
(-1)^{|u|}:u(z)(\Delta v(z)):.
\en
We claim that the following identity holds:
\eq\lb{bracderiv}
\{u(z),:v(z)t(z):\}=:\{u(z),v(z)\}t(z):+(-1)^{(|u|-1)|v|}:v(z)\{u(z),t(z)\}:
\en
ie. the bracket is a derivation in the second argument. This says that $\Delta$
is a second order derivation, and hence proves that it is a BV operator on the
cohomology.
\be{lem}
For any homogeneous quantum operators $A(z),B(z),C(z)$, the following holds:
\eqa\lb{braciden}
A(z)\circ_1(:B(z)C(z):)&-&:(A(z)\circ_1B(z))C(z):
-(-1)^{|A||B|}:B(z)(A(z)\circ_1C(z)):\nnb\\
&=&(A(z)\circ_0B(z))\circ_0C(z).
\ena
\end{lem}
When $A(z)=b(z)$, the LHS of \erf{braciden} is nothing but
$(-1)^{|B|}\{B(z),C(z)\}$ while the RHS is clearly a derivation in the argument
$C(z)$. Thus the lemma implies the identity \erf{bracderiv}.

Proof of lemma: The LHS of \erf{braciden} is:
\eqa
LHS&=&[A(1)-A(0)z,:B(z)C(z):]-:[A(1)-A(0)z,B(z)]C(z):\nnb\\
& &-(-1)^{|A||B|}:B(z)[A(1)-A(0)z,C(z)]:\nnb\\
&=&-[A(0)z,B(z)^-]C(z)-(-1)^{|B||C|+|A||C|}C(z)[A(0)z,B(z)^+]\nnb\\
& &-(-1)^{|A||B|}:B(z)[A(0)z,C(z)]:\nnb\\
& &+[A(0)z,B(z)]^-C(z)+(-1)^{(|A|+|B|)|C|}C(z)[A(0)z,B(z)]^+\nnb\\
& &+(-1)^{|A||B|}:B(z)[A(0)z,C(z)]:\nnb\\
&=&-[A(0)z,B(z)^-]C(z)-(-1)^{|B||C|+|A||C|}C(z)[A(0)z,B(z)^+]\nnb\\
& &+([A(0),B(0)]+[A(0)z,B(z)^-])C(z)\nnb\\
& &+(-1)^{(|A|+|B|)|C|}C(z)([A(0)z,B(z)^+]-[A(0),B(0)])\nnb\\
&=&[A(0),B(0)]C(z)-(-1)^{(|A|+|B|)|C|}C(z)[A(0),B(0)]\nnb\\
&=&[A(0),[B(0),C(z)]]\nnb\\
&=&RHS.
\ena
This proves our lemma and completes our proof of the theorem. $\Box$

The two main theorems above were originally proved in \cite{LZ9} in the
context of vertex operator algebras. (For related
versions of Theorem \ref{5.4}, see \cite{GJ}\cite{SP}\cite{KSV}\cite{Hu}.)

\noi{\footnotesize DEPARTMENT OF MATHEMATICS, HARVARD UNIVERSITY
CAMBRIDGE, MA 02138. lian$@$math.harvard.edu}\\

\noi{\footnotesize DEPARTMENT OF MATHEMATICS, YALE UNIVERSITY
NEW HAVEN, CT 06520. gregg$@$math.yale.edu}

\end{document}